\documentclass[12pt]{article}

\usepackage{latexsym, amssymb, amsmath, amscd, amsfonts, epsfig, graphicx, colordvi,verbatim,ifpdf}
\usepackage{amsfonts, amsmath, amssymb}
\usepackage{amssymb,amsfonts,amsmath,latexsym,epsfig,cite, psfrag,eepic,color}
\usepackage{amscd,graphics}
\usepackage{latexsym, amssymb,  amsmath,amscd, amsfonts, epsfig, graphicx, colordvi,amsthm}
\usepackage{float}

\usepackage{graphicx}
\usepackage{color}
\usepackage{ifpdf}
\usepackage{fancybox}
\usepackage{makecell}
\usepackage{multicol}
\usepackage{multirow}

\newtheorem{thm}{Theorem}[section]

\newtheorem{defi}[thm]{Definition}

\def\pf{\noindent{\it Proof.} }
\def\qed{\nopagebreak\hfill{\rule{4pt}{7pt}}
\medbreak}

\numberwithin{equation}{section}

\def\qed{\nopagebreak\hfill{\rule{4pt}{7pt}}
\medbreak}

\newlength{\boxedparwidth}
  {\begin{center} \begin{tabular}{|@{\hspace{.315in}}c@{\hspace{.15in}}|}
                  \hline \\ \begin{minipage}[t]{\boxedparwidth}
                  \setlength{\parindent}{.25in}}%
  {\end{minipage} \\ \\ \hline \end{tabular} \end{center}}

\allowdisplaybreaks
\begin{document}
\begin{center}

{ \large\bf  A  rank of  partitions with overline designated summands}
\end{center}

\vskip 5mm

\begin{center}
{  Robert X.J. Hao}$^{1}$ ,    {Erin  Y.Y. Shen}$^{2}$
, {Wenston J.T. Zang}$^{3}$  \vskip 2mm

   $^{1}$ Department of Mathematics and Physics, \\Nanjing Institute Of Technology,
    Nanjing 211167, P.~R.~China\\[6pt]
    $^{2}$School of Science, \\ Hohai University,
    Nanjing 210098, P.~R.~China\\[6pt]
   $^{3}$Institute for Advanced Study in Mathematics,\\
Harbin Institute of Technology, Heilongjiang 150001, P.~R.~China\\[6pt]

   \vskip 2mm

    $^1$haoxj@njit.edu.cn, $^2$shenyy@hhu.edu.cn, $^3$zang@hit.edu.cn
\end{center}

\vskip 6mm \noindent {\bf Abstract.}
Andrews, Lewis and Lovejoy introduced the
partition function $PD(n)$ as the number of partitions of $n$ with designated summands. In a recent work, Lin studied a partition function $PD_{t}(n)$ which counts the number of tagged parts over all the partitions of $n$ with designated summands. He proved that  $PD_{t}(3n+2)$ is divisible by $3$. In this paper, we first introduce a structure named partitions with overline designated summands, which is counted by $PD_t(n)$. We then define a  generalized rank  of partitions with overline designated summands and give a combinatorial interpretation of the congruence for $PD_t(3n+2)$.



\section{Introduction}

 In \cite{Andrews-Lewis-Lovejoy-2002}, Andrews, Lewis and Lovejoy
defined partitions with designated summands on ordinary partitions by designating exactly one part among parts with equal size. They studied the number of partitions with designated summands.
 Let $PD(n)$ denote the number of partitions of $n$ with designated summands.
  For example, there are 15 partitions of $5$ with designated summands:
{\footnotesize\begin{align*}
\begin{array}{ccccc}
5',      & 4'+1',  & 3'+2',    & 3'+1'+1,    & 3'+1+1',  \\[5pt]
2'+2+1', & 2+2'+1',  & 2'+1'+1+1,& 2'+1+1'+1,& 2'+1+1+1', \\[5pt]
1'+1+1+1+1, & 1+1'+1+1+1, & 1+1+1'+1+1, & 1+1+1+1'+1, & 1+1+1+1+1'.
\end{array}
\end{align*}}
Thus $PD(5)=15$. Andrews, Lewis and Lovejoy\cite{Andrews-Lewis-Lovejoy-2002}
obtained the generating function of $PD(n)$ as given by
\begin{align}\label{PD(n)}
\sum_{n=0}^{\infty}PD(n)q^n=\frac{(q^6;q^6)_{\infty}}{(q;q)_{\infty}
(q^2;q^2)_{\infty}(q^3;q^3)_{\infty}},
\end{align}
where here and throughout this paper,
$(a;q)_\infty$ stands for the $q$-shifted factorial
\[
(a;q)_{\infty}=\prod_{n=1}^{\infty}(1-aq^{n-1}),\,\, |q|<1.
\]
By using modular forms and $q$-series identities,  Andrews, Lewis and Lovejoy\cite{Andrews-Lewis-Lovejoy-2002} proved some  arithmetic properties of the partition function $PD(n)$. For instance, they obtained
a Ramanujan-type congruence as given by
\begin{align}
\label{modPD(3n+2)}PD(3n+2)\equiv 0\pmod{3}.
\end{align}
By introducing the $pd$-rank for  partitions with designated summands, Chen, Ji, Jin and the second author\cite{Chen-Ji-Jin-Shen-2013} provided  a combinatorial interpretation for the congruence \eqref{modPD(3n+2)}.

Recently, Lin\cite{Lin-2018} introduced a partition function $PD_t(n)$, which counts the number of tagged parts over all the partitions of
$n$ with designated summands. There are 24 tagged parts over all the partitions of $5$ with designated summands. Hence $PD_t(5)=24$. Lin\cite{Lin-2018} showed that the generating function for $PD_t(n)$ is
\begin{align}\label{PDt(n)-1}
\sum_{n=0}^{\infty}PD_t(n)q^n&=\frac{(q^6;q^6)_{\infty}}{(q;q)_{\infty}
(q^2;q^2)_{\infty}(q^3;q^3)_{\infty}}\sum_{k\geq 1}\frac{q^k+q^{2k}}{1+q^{3k}}\\[5pt]\label{PDt(n)-2}
&=\frac{(q^6;q^6)_{\infty}}{(q;q)_{\infty}
(q^2;q^2)_{\infty}(q^3;q^3)_{\infty}}
\left(\frac{(q^3;q^3)^6_{\infty}(q^2;q^2)_{\infty}}
{2(q^6;q^6)^3_{\infty}(q;q)^2_{\infty}}-\frac{1}{2}\right).
\end{align}
Lin\cite{Lin-2018} also established many congruences modulo small powers
of $3$ for the partition function $PD_t(n)$ including a Ramanujan-type
congruence as given by
\begin{align}
\label{modPDt(3n+2)}PD_t(3n+2)\equiv0\pmod{3}.
\end{align}
In the end of his paper, Lin asked for a suitable rank
of partitions with designated summands that could combinatorially
interpret the congruence \eqref{modPDt(3n+2)}. This is our main task.

Recall that for a partition $\lambda=(\lambda_1,\lambda_2,\ldots,\lambda_\ell)$, the crank of $\lambda$     is defined  by
 \[\text{crank}(\lambda)=\left\{
\begin{array}{ll}
\lambda_1,\ \ &\text{ if }\ n_1(\lambda)=0,\\[10pt]
 \mu(\lambda)-n_1(\lambda), \ &\text{ if }\ n_1(\lambda)>0,
\end{array}\right.
\]
where   $n_1(\lambda)$ is  the number of ones in $\lambda$ and $\mu(\lambda)$ is  the number of parts  larger than $n_1(\lambda)$.
Let $M(m,n)$ denote the number of partitions of $n$ with crank $m$. We use the convention that
\begin{equation}\label{equ-mn-1}
M(1,1)=M(-1,1)=1,\quad M(0,1)=-1,\quad M(i,1)=0,
\end{equation}
for all $i\neq 0,1,-1$.
Andrews and Garvan \cite{Andrews-Garvan-1988} gave the following generating function of $M(m,n)$ as
\begin{equation}\label{equ-gef-mmn}
\sum_{m=-\infty}^\infty\sum_{n=0}^\infty M(m,n)z^m q^n=\frac{(q;q)_\infty}{(zq;q)_\infty(z^{-1}q;q)_\infty}.
\end{equation}
In this paper, we first introduce partitions with overline designated summands which are counted by the partition function $PD_t(n)$. Then we define a statistic called the $pdt$-rank of partitions with overline designated summands.
Let $N_{dt}(m,n)$ denote the number of  partitions of $n$ with overline designated summands with $pdt$-rank $m$ and let $$N_{dt}(i,3;n)=\sum_{m\equiv i\pmod{3}} N_{dt}(m,n).$$  We obtain the following theorem. It provides a combinatorial interpretation for the congruence \eqref{modPDt(3n+2)}.

\begin{thm}\label{thm-c}
For $n\geq 0$,  we have
\begin{align}
N_{dt}(0,3;3n+2)=N_{dt}(1,3;3n+2)=N_{dt}(2,3;3n+2)=\frac{PD_t(3n+2)}{3}.
\end{align}
\end{thm}

This paper is organized as follows. First we give the definition of the $pdt$-rank of partitions with overline designated summands in Section $2$. By using the
$pdt$-rank, we provide a proof for Theorem \ref{thm-c} in Section $3$.
In Section $4$, we introduce the modified $pdt$-rank which enables us to
divide the set of partitions with  overline designated summands counted by $PD_t(3n+2)$ into three equinumerous subsets.

\section{The $pdt$-rank}

In this section, we aim to introduce a rank called the $pdt$-rank.
In doing so, we first define  partitions with overline designated summands. A partition with overline designated summands is defined on a partition with  designated summands with exactly one designated part be overlined. It is clear that the number of partitions with overline designated summands  of $n$ is equal to $PD_t(n)$. For instance, there are $24$  partitions of $5$ with overline designated summands:
{\footnotesize\begin{align*}
\begin{array}{ccccc}
\overline{5'},      & \overline{4'}+1', & 4'+\overline{1'}, & \overline{3'}+2',   & 3'+\overline{2'},  \\[5pt]
\overline{3'}+1'+1,  & 3'+\overline{1'}+1,    & \overline{3'}+1+1', & 3'+1+\overline{1'}, & \overline{2'}+2+1',\\[5pt]
2'+2+\overline{1'}, & 2+\overline{2'}+1', & 2+2'+\overline{1'}, & \overline{2'}+1'+1+1,& 2'+\overline{1'}+1+1,\\[5pt]
\overline{2'}+1+1'+1, & 2'+1+\overline{1'}+1,& \overline{2'}+1+1+1', & 2'+1+1+\overline{1'}, & \overline{1'}+1+1+1+1,\\[5pt]
1+\overline{1'}+1+1+1, & 1+1+\overline{1'}+1+1, & 1+1+1+\overline{1'}+1, & 1+1+1+1+\overline{1'}.
\end{array}
\end{align*}}
Therefore $PD_t(5)=24$.

Here we  introduce a different way to denote partitions with overline designated summands. It is well known that a partition $\lambda$ of $n$ can be written as ($1^{f_1}2^{f_2}\ldots n^{f_n}$), where $f_i\geq 0$ denotes the number of $i$'s appears in $\lambda$. For a partition $\lambda$ with overline designated summands, we denote $\lambda$ as $(1^{f_1}2^{f_2}\ldots n^{f_n},g_1,g_2,\ldots,g_n;k)$. For $f_i=0$, we set $g_i=0$.
For $f_i\geq 1$, the $g_i$th $i$ is designated (from left to right), and the designated part $k'$ is overlined. Hence when $f_i \geq 1$, we have $1\leq g_i\leq f_i$ and $f_k\ge 1$.
 For instance, the partition with overline designated summands $\overline{3'}+1+1'$ can be written as $ (1^23^1,2,0,1;3)$.

In order to define the $pdt$-rank of partitions with overline designated summands,
we proceed to  build a bijection  between the set of partitions of $n$ with overline designated summands and a set of pairs of partitions as follows.

Let $S_1(n)$ denote the set of partitions with overline designated summands $(1^{f_1}2^{f_2}\ldots n^{f_n},\break g_1,g_2,\ldots,g_n;k)$ such that $\sum_{i=1}^n if_i=n$. Let $S_2(n)$ denote the set of triplets $(\alpha,\beta;t)$, where $\alpha=(1^{x_1}2^{x_2}\ldots n^{x_n})$ and $\beta=(1^{y_1}2^{y_2}\ldots n^{y_n})$ are ordinary partitions, and $t$ is a positive integer such that $x_t\ge 1$ and $x_i\neq 1$ for all $i\neq t$. Moreover, $\sum_{i=1}^n i(x_i+y_i)=n$. We have the following result.

\begin{thm}\label{thm-delta}
There  is a bijection $\Delta$ between $S_1(n)$ and $S_2(n)$.
\end{thm}

\pf  Let $\lambda=(1^{f_1}2^{f_2}\ldots n^{f_n},g_1,g_2,\ldots,g_n;k)\in S_1(n)$ be a partition of $n$ with overline designated summands. For each $1\leq i\leq n$, we define nonnegative integers $x_i$ and $y_i$ as given below. There are three cases.
\begin{itemize}
\item[Case 1.] If $f_i=0$, then set $x_i=y_i=0$;
\item[Case 2.] If $i=k$ or $f_i\geq g_i\geq 2$, then set $x_i=g_i$ and $y_i=f_i-g_i$;
\item[Case 3.] If $i\neq k$ and $g_i=1$, then set $x_i=0$ and $y_i=f_i$.
\end{itemize}

Let $\alpha=(1^{x_1}2^{x_2}\ldots n^{x_n})$ and $ \beta= (1^{y_1}2^{y_2}\ldots n^{y_n})$. It is easy to check that $x_i+y_i=f_i$ for $1\leq i\leq n$ and $x_k\geq 1$ and $x_i\neq 1$ for all $i\neq k$. Hence we may set $\Delta(\lambda)=(\alpha,\beta;k)\in S_2(n)$.

To show that $\Delta$ is a bijection, we need to construct the inverse map $\Delta^{-1}$. Given  $(\alpha,\beta;t)\in S_2(n)$, where $\alpha=(1^{x_1}\ldots n^{x_n})$ and $\beta=(1^{y_1}\ldots n^{y_n})$, we define $f_i$ and $g_i$ with $1\leq i\leq n$ as follows.
\begin{itemize}
\item[Case 1.] If $x_i=y_i=0$, then set $f_i=g_i=0$;
\item[Case 2.] If $i=t$ or $x_i\geq 2$, then set $f_i=x_i+y_i$ and $g_i=x_i$;
\item[Case 3.] If $x_i=0$ and $y_i\geq 1$, then set $f_i=x_i+y_i$ and $g_i=1$.
\end{itemize}

Let $\lambda=(1^{f_1}2^{f_2}\ldots n^{f_n},g_1,g_2,\ldots,g_n;t)$,  it can be checked that if $f_i=0$ we get $g_i=0$ according to the above construction. Furthermore, $f_t\ge 1$ and when $f_i\ge 1$, we have $f_i\ge g_i\ge 1$. Moreover $f_i=x_i+y_i$, so we have $\sum_{i=1}^n if_i=\sum_{i=1}^n i(x_i+y_i)$. Thus $\lambda\in S_1(n)$.

 It is clear to see that $\Delta^{-1}$ is the inverse map of $\Delta$. Thus $\Delta$ is a bijection. This completes the proof.\qed

 For instance, let $\lambda=(5',5,3,3',3,\overline{2}',2,1')$, then we may denote $\lambda$ as $(1^1,2^2,3^3,5^2,1,1,2,$
 $0,1;2)$. Applying the bijection  $\Delta$ to $\lambda$, we derive that $\alpha=(2^1,3^2)$ and $\beta=(1^1,2^1,3^1,5^2)$. Thus $\Delta(\lambda)=((2^1,3^2),(1^1,2^1,3^1,5^2);2)$. Applying $\Delta^{-1}$ to $((2^1,3^2),(1^1,2^1,3^1,5^2);2)$, we recover $\lambda$.

We are now in a position to define the $pdt$-rank.

\begin{defi}
Let $\lambda$ be a partition with  overline designated summands and let $\Delta(\lambda)=(\alpha,\beta;k)$. The $pdt$-rank of $\lambda$, denoted  $r_{dt}(\lambda)$, is defined by
\begin{align}
r_{dt}(\lambda)=crank(\beta),
\end{align}
where $crank(\beta)$ is the crank of partition $\beta$.
\end{defi}

Recall that $N_{dt}(m,n)$ denote the number of partitions of $n$ with overline designated summands with $pdt$-rank $m$. Here we make the appropriate modifications based on the fact that for ordinary
partitions  ${M}(0,1)=-1$ and ${M}(-1,1) ={M}(1,1) = 1$. For example, the partition with overline designated summands ${2}+\overline{2'}+1'$
can be divided into $\alpha=(2,2)$ and $\beta=(1)$ under the bijection $\Delta$. When $\beta=(1)$, we use the convention that this partition contributes a $-1$ to the count of $N_{dt}(0, 5)$, a $1$ to $N_{dt}(-1, 5)$ and $N_{dt}(1, 5)$ respectively.

For example, Table \ref{table-1} gives the 24 partitions of $5$ with overline designated summands. According to  this example, we know that four overline designated partitions, namely,   $\bar{4}'+1'$, $2+\bar{2}'+1$, $\bar{2}'+1+1'+1$ and $1+1+1+\bar{1}'+1$ are divided into $\alpha,\beta$
with $\beta=(1)$. By the definition of $pdt$-rank, we have $N_{dt}(0, 5)=8-4=4$, $N_{dt}(-1, 5)=0+4=4$ and $N_{dt}(1, 5)=0+4=4$.

{\makeatletter\def\@captype{table}\makeatother\centering
\begin{tabular}{m{4cm}<{\centering}|m{4cm}<{\centering}|m{1.7cm}<{\centering}}
$\lambda$       & $(\alpha,\beta)$    &$r_{dt}(\lambda) $\\[4pt]
\hline
$\bar{5}'$      & $(5,\ \emptyset)$
&$0$\\
\hline
$\bar{4}'+1'$             &  $(4,\ 1)$       & $-$\\
\hline
$4'+\bar{1}'$             &  $(1,\ 4)$       & $4$\\
\hline
$\bar{3}'+2'$             &  $(3,\ 2)$       & $2$\\
\hline
$3'+\bar{2}'$             &  $(2,\ 3)$       & $3$\\
\hline
$\bar{3}'+1'+1$      &  $(3,\ 1+1)$       & $-2$\\
\hline
$3'+\bar{1}'+1$      &  $(1,\ 3+1)$       & $0$\\
\hline
$\bar{3}'+1+1'$     &  $(3+1+1,\ \emptyset)$       & $0$\\
\hline
$3'+1+\bar{1}'$     &  $(1+1,\ 3)$       & $3$\\
\hline
$\bar{2}'+2+1'$     &  $(2,\ 2+1)$       & $0$\\
\hline
$2'+2+\bar{1}'$     &  $(1,\ 2+2)$       & $2$\\
\hline
$2+\bar{2}'+1'$     &  $(2+2,\ 1)$       & $-$\\
\hline
$2+2'+\bar{1}'$     &  $(2+2+1,\emptyset )$       & $0$\\
\hline
$\bar{2}'+1'+1+1$     &  $(2,\ 1+1+1)$       & $-3$\\
\hline
$2'+\bar{1}'+1+1$     &  $(1,\ 2+1+1)$       & $-2$\\
\hline
$\bar{2}'+1+1'+1$     &  $(2+1+1,\ 1)$       & $-$\\
\hline
$\bar{2}'+1+1+1'$     &  $(2+1+1+1,\ \emptyset)$       & $0$\\
\hline
$2'+1+1+\bar{1}'$     &  $(1+1+1,\ 2)$       & $2$\\
\hline
$2'+1+\bar{1}'+1$     &  $(1+1,\ 2+1)$       & $0$\\
\hline
$\bar{1}'+1+1+1+1$     &  $(1,\ 1+1+1+1)$       & $-4$\\
\hline
$1+\bar{1}'+1+1+1$     &  $(1+1,\ 1+1+1)$       & $-3$\\
\hline
$1+1+\bar{1}'+1+1$     &  $(1+1+1,\ 1+1)$       & $-2$\\
\hline
$1+1+1+\bar{1}'+1$     &  $(1+1+1+1,\ 1)$       & $-$\\
\hline
$1+1+1+1+\bar{1}'$     &  $(1+1+1+1+1,\ \emptyset)$       & $0$\\
\end{tabular}\caption{ The case for $n=5$ with $pdt$-rank $r_{dt}(\lambda)$.}\label{table-1}
}

We next derive the generating function of $N_{dt}(m,n)$. By the definition, it is clear that the generating function of $\alpha$ equals
\begin{equation}
\sum_{k=1}^{\infty}
\frac{q^k}{1-q^k}\prod_{j=1,j\not=k}^{\infty}
\left(1+\frac{q^{2j}}{1-q^{j}}\right).
\end{equation}
Since $\beta$ is an ordinary partition and the $pdt$-rank only relies on $\beta$, by \eqref{equ-gef-mmn}, the generating function of $N_{dt}(m,n)$ can be given as
\begin{equation}\label{equ-gef=ndtmn}
\sum_{n=0}^{\infty}\sum_{m=-\infty}^{\infty}N_{dt}(m,n)z^mq^n
=\frac{(q;q)_{\infty}}{(zq;q)_{\infty}(z^{-1}q;q)_{\infty}}
\sum_{k=1}^{\infty}
\frac{q^k}{1-q^k}\prod_{j=1,j\not=k}^{\infty}
\left(1+\frac{q^{2j}}{1-q^{j}}\right).
\end{equation}

\section{A proof of Theorem \ref{thm-c}}

In this section, we provide a proof for Theorem \ref{thm-c}.

\noindent{\it Proof of Theorem \ref{thm-c}.}  By \eqref{equ-gef=ndtmn}, we have
\begin{align}
\sum_{n=0}^{\infty}\sum_{m=-\infty}^{\infty}N_{dt}(m,n)z^mq^n
&=\frac{(q;q)_{\infty}}{(zq;q)_{\infty}(z^{-1}q;q)_{\infty}}
\sum_{k=1}^{\infty}
\frac{q^k}{1-q^k}\prod_{j=1,j\not=k}^{\infty}
\left(1+\frac{q^{2j}}{1-q^{j}}\right)\nonumber \\
&=\frac{(q;q)_{\infty}}{(zq;q)_{\infty}(z^{-1}q;q)_{\infty}}
\sum_{k=1}^{\infty}\frac{q^k}{1-q^k}\prod_{j=1,j\not=k}^{\infty}
\frac{1+q^{3j}}{1-q^{2j}}
\nonumber \\
&=\frac{(q;q)_{\infty}}{(zq;q)_{\infty}(z^{-1}q;q)_{\infty}}
\sum_{k=1}^{\infty}\frac{q^k}{1-q^k}\frac{(-q^3;q^3)_\infty}{(q^2;q^2)_\infty}
\frac{1-q^{2k}}{1+q^{3k}}
\nonumber \\
&=\frac{(q;q)_{\infty}}{(zq;q)_{\infty}(z^{-1}q;q)_{\infty}}
\frac{(-q^3;q^3)_{\infty}}{(q^2;q^2)_{\infty}}
\sum_{k=1}^{\infty}
\frac{q^k+q^{2k}}{1+q^{3k}}
\nonumber \\
&=\frac{(q;q)_{\infty}}{(zq;q)_{\infty}(z^{-1}q;q)_{\infty}}
\frac{(q^6;q^6)_{\infty}}{(q^2;q^2)_{\infty}(q^3;q^3)_{\infty}}
\sum_{k=1}^{\infty}\frac{q^k+q^{2k}}{1+q^{3k}}.\label{temp-pdt}
\end{align}
Recall that Lin \cite[Theorem 3.2]{Lin-2018} derive the following identity
\begin{equation}\label{equ-Lin-32}
\sum_{k=1}^{\infty}\frac{q^k+q^{2k}}{1+q^{3k}}
=\sum_{k=-\infty}^{\infty}\frac{q^k}{1+q^{3k}}-\frac{1}{2}
=\frac{(q^3;q^3)^6_{\infty}(q^2;q^2)_{\infty}}
{2(q^6;q^6)^3_{\infty}(q;q)^2_{\infty}}-\frac{1}{2}.
\end{equation}
Setting $z=\zeta=e^{\frac{2\pi i}{3}}$ in \eqref{temp-pdt}, we derive that
\begin{align*}
\sum_{n=0}^{\infty}
\sum_{m=-\infty}^{\infty}N_{dt}(m,n)\zeta^mq^n
&=\sum_{n=0}^{\infty}\sum_{i=0}^{2}N_{dt}(i,3;n)\zeta^iq^n
\nonumber \\
&=\frac{(q;q)_{\infty}}{(\zeta q;q)_{\infty}(\zeta^{-1}q;q)_{\infty}}
\frac{(q^6;q^6)_{\infty}}{(q^2;q^2)_{\infty}(q^3;q^3)_{\infty}}
\sum_{k=1}^{\infty}\frac{q^k+q^{2k}}{1+q^{3k}}.
\end{align*}
Using \eqref{equ-Lin-32}, we find that
\begin{equation}\label{GF3}
\sum_{n=0}^{\infty}
\sum_{m=-\infty}^{\infty}N_{dt}(m,n)\zeta^mq^n=\frac{(q;q)_{\infty}}{(\zeta q;q)_{\infty}(\zeta^{-1}q;q)_{\infty}}
\frac{(q^6;q^6)_{\infty}}{(q^2;q^2)_{\infty}(q^3;q^3)_{\infty}}
\left(\frac{(q^3;q^3)_{\infty}^6(q^2;q^2)_{\infty}}
{2(q^6;q^6)_{\infty}^3(q;q)_{\infty}^2}-\frac{1}{2}\right).
\end{equation}
Multiplying the right side of \eqref{GF3} by
\begin{equation*}
\frac{(q;q)_{\infty}}{(q;q)_{\infty}},
\end{equation*}
and noting that
\begin{equation*}
(1-x)(1-x\zeta)(1-x\zeta^2)=1-x^3,
\end{equation*}
we derive that
\begin{align}\label{PDD012}
\sum_{n=0}^{\infty}\sum_{i=0}^{2}N_{dt}(i,3;n)\zeta^iq^n
&=\frac{(q;q)_{\infty}^2}{(\zeta q;q)_{\infty}
(\zeta^{-1}q;q)_{\infty}(q;q)_{\infty}}
\frac{(q^6;q^6)_{\infty}}{(q^2;q^2)_{\infty}(q^3;q^3)_{\infty}}
\left(\frac{(q^3;q^3)_{\infty}^6(q^2;q^2)_{\infty}}
{2(q^6;q^6)_{\infty}^3(q;q)_{\infty}^2}-\frac{1}{2}\right)\nonumber \\[5pt]
&=\frac{1}{2}\left(\frac{(q^3;q^3)_{\infty}^4}{(q^6;q^6)_{\infty}^2}-
\frac{1}{(q^3;q^6)_{\infty}^2(q^6;q^6)_{\infty}}\cdot
\frac{(q;q)_{\infty}^2}{(q^2;q^2)_{\infty}}\right).
\end{align}
Using Jacobi triple product identity \cite[Theorem 1.3.3]{Berndt-2006}
\[
\sum_{n=-\infty}^{\infty}z^nq^{n^2}
=(-zq;q^2)_{\infty}(-q/z;q^2)_{\infty}(q^2;q^2)_{\infty}
\]
with $z=-1$, we have
\begin{equation}\label{eq-frac-q-q-infty}
\frac{(q;q)_{\infty}^2}{(q^2;q^2)_{\infty}}
=(q;q^2)_{\infty}^2(q^2;q^2)_{\infty}
=\sum_{n=-\infty}^{\infty}(-1)^nq^{n^2},
\end{equation}
Substituting \eqref{eq-frac-q-q-infty} into \eqref{PDD012}, we obtain that
\begin{align}\label{pdt-nn}
\sum_{n=-\infty}^{\infty}\sum_{i=0}^2N_{dt}(i,3;n)\zeta^i q^n=\frac{1}{2}\left(\frac{(q^3;q^3)_{\infty}^4}{(q^6;q^6)_{\infty}^2}-
\frac{1}{(q^3;q^6)_{\infty}^2(q^6;q^6)_{\infty}}
\sum_{n=-\infty}^{\infty}(-1)^nq^{n^2}\right).
\end{align}
Since
\[n^2 \equiv 0 \  \text{ or } \    1 \pmod 3,\]
       the coefficient of
$q^{3n+2}$ in   \eqref{pdt-nn} is zero.
It follows that
\[N_{dt}(0,3;3n+2)+N_{dt}(1,3;3n+2)\zeta+N_{dt}(1,3;3n+2)\zeta^2=0. \]
Since the minimal polynomial of $\zeta$ is $1+x+x^2$,  we conclude that
\[
N_{dt}(0,3;3n+2)=N_{dt}(1,3;3n+2)=N_{dt}(2,3;3n+2).
\] This completes the proof.\qed

\section{The Modified  $pdt$-rank}

Recall that the bijection $\Delta(\lambda)=(\alpha,\beta;k)$,  when $\beta=(1)$, the $pdt$-rank of $\lambda$ contributes a $-1$ to $N_{dt}(0,n)$ and contributes a $1$ to $N_{dt}(-1,n)$ and $N_{dt}(1,n)$ respectively. Thus the $pdt$-rank cannot divide the set of partitions of $3n+2$ with overline designated summands into three equinumerous subsets. In this section, we shall define a modified $pdt$-rank such that the function $N_{dt}(m,n)$ directly counts the number of partitions of $n$ with overline designated summands with modified $pdt$-rank $m$. This statistic enables us to divide the set of partitions of $3n+2$ with overline designated summands into three equinumerous subsets.

In order to define the modified $pdt$-rank, we first define two set of partitions with overline designated summands, namely $A(n)$ and $B(n)$. Here $B(n)$ is the set of partitions $\lambda$ with overline designated summands of $n$ such that $\Delta(\lambda)=(\alpha,(1);k)$. $A(n)$ is a subset of the set of partitions with overline designated summands of $n$ whose $pdt$-rank equals $0$, which will be defined later.
We shall build a bijection $\phi$ between $A(n)$ and $B(n)$, which implies $\# A(n)=\# B(n)$.
After that, we can define the modified $pdt$-rank $r_{mdt}(\lambda)$ as Definition
\ref{rmdtDef}.
Let $N_{mdt}(m,n)$ denote the number of partitions $\lambda$ of $n$ with overline designated summands satisfies that $r_{mdt}(\lambda)=m$. From the above construction, we may see that
\begin{align*}
N_{mdt}(0,n)&=\#\{|\lambda|=n\colon r_{dt}(\lambda)=0\}-\# A(n);\\
N_{mdt}(1,n)&=\#\{|\lambda|=n\colon r_{dt}(\lambda)=1\}+\# A(n);\\
N_{mdt}(-1,n)&=\#\{|\lambda|=n\colon r_{dt}(\lambda)=-1\}+\# B(n).
\end{align*}
Hence we have  $N_{mdt}(m,n)=N_{dt}(m,n)$ for all $m,n$.

We proceed to define the set $A(n)$. To this end,
 we need to  define five sets $A_1(n)$, $A_2(n)$, $A_3(n)$, $A_4(n)$ and  $A_5(n)$ which satisfy  for any $\lambda\in A_i(n)$, $1\leq i\leq 5$, $r_{dt}(\lambda)=0$.  The set $A_1(n)$ is defined as follows.

\begin{defi}\label{DefA1}
Let $A_1(n)$ be the set of partitions of $n$ with overline designated summands
 with the following restrictions:
\begin{itemize}
\item[(1)] $f_i=g_i\neq 1$ for all $i$;
\item[(2)] $k\neq 1$ and $f_1\geq 3$.
\end{itemize}
\end{defi}

We next give the definition of the set $A_2(n)$.

\begin{defi}\label{DefA2}
Let $A_2(n)$ be the set of partitions of $n$ with overline designated summands with the following restrictions:
\begin{itemize}
\item[(1)] $f_i=g_i\neq 1$ for all $i$;
\item[(2)] $k=1$ and $f_1\geq 2$.
\end{itemize}
\end{defi}

The set $A_3(n)$ can be defined as given below.

\begin{defi}\label{DefA3}
Let $A_3(n)$ be the set of partitions of $n$ with overline designated summands with the following restrictions:
\begin{itemize}
\item[(1)] $f_1=g_1=1$;
\item[(2)] $k\neq 1$, $f_k\geq 2$ and $g_k=f_k-1$;
\item[(3)] For all $i\neq k$, $f_i=g_i\neq 1$.
\end{itemize}
\end{defi}

We proceed to define the set $A_4(n)$.

\begin{defi}\label{DefA4}
Let $A_4(n)$ be the set of partitions of $n$ with overline designated summands with the following restrictions:
\begin{itemize}
\item[(1)] $f_1=g_1=1$;
\item[(2)] $k\neq 1$ and $f_k=g_k \geq 2$;
\item[(3)] For all $i$, $f_i=g_i$ and there exists a unique $j\neq 1$ such that $f_j=g_j=1$.
\end{itemize}
\end{defi}

Finally, we define the set $A_5(n)$.

\begin{defi}\label{DefA5}
Let $A_5(n)$ be the set of partitions of $n$ with overline designated summands with the following restrictions:
\begin{itemize}
\item[(1)] $k\neq 1$ and $f_k=g_k=1$;
\item[(2)] For all $i\neq k$, $f_i=0$.
\end{itemize}
\end{defi}

It is trivial to check that for any $\lambda\in A_i(n)$, $1\leq i\leq 5$, $r_{dt}(\lambda)=0$ and  $A_i(n)$ are disjoint.
 Let $A(n)=A_1(n)\cup A_2(n)\cup A_3(n)\cup A_4(n)\cup A_5(n)$. Clearly,  $A(n)$ is a subset of the set of partitions  $\lambda$ of $n$ with overline designated summands such that $r_{dt}(\lambda)=0$.

To establish a bijection $\phi$ between $A(n)$ and $B(n)$, we divide
$B(n)$ into five disjoint subsets $B_i$ for $1\leq i\leq 5$. We then construct five bijections $\phi_i$ between $A_i(n)$ and $B_i(n)$, $1\leq i\leq 5$. We now give the definitions of $B_i(n)$ for $1\leq i\leq 5$.

\begin{defi}\label{DefB1}
Let $B_1(n)$ be the set of partitions of $n$ with overline designated summands with the following restrictions:
\begin{itemize}
\item[(1)] $f_i=g_i$ for all $i\neq 1$;
\item[(2)] $k\neq 1$, $f_1\geq 3$ and $g_1=f_1-1$;
\item[(3)] $f_i\neq 1$ for all $i\neq 1,k$.
\end{itemize}
\end{defi}
\begin{defi}\label{DefB2}
Let $B_2(n)$ be the set of partitions of $n$ with overline designated summands with the following restrictions:
\begin{itemize}
\item[(1)] $f_i=g_i\neq 1$ for all $i\neq 1$;
\item[(2)] $k=1$, $f_1\geq 2$ and $g_1=f_1-1$.
\end{itemize}
\end{defi}
\begin{defi}\label{DefB3}
Let $B_3(n)$ be the set of partitions of $n$ with overline designated summands with the following restrictions:
\begin{itemize}
\item[(1)] $f_1=g_1=1$ and $f_i=g_i\neq 1$ for all $i\neq 1$;
\item[(2)] $k\neq 1$.
\end{itemize}
\end{defi}

\begin{defi}\label{DefB4}
Let $B_4(n)$ be the set of partitions of $n$ with overline designated summands with the following restrictions:
\begin{itemize}
\item[(1)] $f_1=g_1=1$;
\item[(2)] $k\neq 1$ and $f_k=g_k=1$;
\item[(3)] For all $i\neq k$, $f_i=g_i\neq 1$ and there exists $j\neq 1,k$ such that $f_j=g_j\geq 2$.
\end{itemize}
\end{defi}

\begin{defi}\label{DefB5}
Let $B_5(n)$ be the set of partitions of $n$ with overline designated summands with the following restrictions:
\begin{itemize}
\item[(1)] $k\neq 1$ and $f_k=g_k=1$;
\item[(2)] $f_1=g_1=1$;
\item[(3)] For all $i\neq 1,k$, $f_i=0$.
\end{itemize}
\end{defi}

It can be checked that
\[B(n)=\biguplus_{i=1}^5 B_i(n).\]
We are now in a position to  present the five bijections $\phi_i$ between $A_i(n)$ and $B_i(n)$ for $1\leq i\leq 5$.

\begin{thm}\label{thn=lem-1}
There is a bijection $\phi_1$ between  $A_1(n)$ and $B_1(n)$.
\end{thm}
\pf For any $\lambda_1=(1^{f_1}2^{f_2}\ldots n^{f_n},g_1,g_2,\ldots,g_n;k)\in A_1(n)$, by Definition \ref{DefA1}, we have $f_i=g_i\neq 1$, $k\neq 1$ and $f_1\geq 3$.
Let $\phi_1(\lambda_1)=\mu_1=(1^{f_1}2^{f_2}\ldots n^{f_n},g_1-1,g_2,\ldots,g_n;k)$. It is clear that $\mu_1\in B_1(n)$ and $\phi_1$ is a bijection. This completes the proof.\qed

\begin{thm}\label{thn=lem-2}
There is a bijection $\phi_2$ between  $A_2(n)$ and $B_2(n)$.
\end{thm}
\pf For any $\lambda_2=(1^{f_1}2^{f_2}\ldots n^{f_n},g_1,g_2,\ldots,g_n;k)\in A_2(n)$, by Definition \ref{DefA2}, we have $f_i=g_i\neq 1$, $k=1$ and $f_1\geq 2$.
Let $\phi_2(\lambda_2)=\mu_2=(1^{f_1}2^{f_2}\ldots n^{f_n},g_1-1,g_2,\ldots,g_n;k)$. It is easy to check that $\mu_2\in B_2(n)$ and $\phi_2$ is a bijection. This completes the proof.\qed

\begin{thm}\label{thn=lem-3}
There is a bijection $\phi_3$ between  $A_3(n)$ and $B_3(n)$.
\end{thm}
\pf For any $\lambda_3=(1^{f_1}2^{f_2}\ldots n^{f_n},g_1,g_2,\ldots,g_n;k)\in A_3(n)$, by Definition \ref{DefA3}, we have $f_1=g_1=1$, $f_k\geq 2$, $g_k=f_k-1$, and $f_i=g_i\neq 1$ for all $i\neq 1,k$.
Let $\phi_3(\lambda_3)=\mu_3=(1^{f_1}2^{f_2}\ldots n^{f_n},g_1,g_2,\ldots,g_k+1,\ldots, g_n;k)$. We see that $g_k+1=f_k$, this implies that $f_i=g_i\neq 1$ for all $i\neq 1$ and $f_1=g_1=1$. Hence $\mu_3\in B_3(n)$ and it is clear that $\phi_3$ is a bijection. This completes the proof.\qed

\begin{thm}\label{thn=lem-4}
There is a bijection $\phi_4$ between  $A_4(n)$ and $B_4(n)$.
\end{thm}
\pf For any $\lambda_4=(1^{f_1}2^{f_2}\ldots n^{f_n},g_1,g_2,\ldots,g_n;k)\in A_4(n)$, by Definition \ref{DefA4}, we have $f_i=g_i$ for all $i$, especially, $f_1=g_1=1$, $g_k=f_k\geq 2$, and there exists a unique $j\neq 1$ such that $f_j=g_j=1$.
Let $\phi_4(\lambda_4)=\mu_4=(1^{f_1}2^{f_2}\ldots n^{f_n},g_1,g_2,\ldots, g_n;j)$. We next check that $\mu_4\in B_4$.
Clearly, $f_1=g_1=f_j=g_j=1$, and for all $i\neq 1,j$, $f_i=g_i\neq 1$. Moreover, $f_k=g_k\geq 2$. Hence, $\mu_4\in B_4$. It is trivial to check that $\phi_4$ is a bijection. This completes the proof.\qed

\begin{thm}\label{thn=lem-5}
There is a bijection $\phi_5$ between  $A_5(n)$ and $B_5(n)$.
\end{thm}
\pf For fixed $n$, it is clear that $A_5(n)$ has only one element $(\overline{n'})$ and $B_5(n)$ only contains one element $(\overline{n-1'},1')$. Set $\phi_5(\overline{n'})=(\overline{n-1'},1')$ and the proof is completed.\qed

Combining Theorem \ref{thn=lem-1}, \ref{thn=lem-2}, \ref{thn=lem-3}, \ref{thn=lem-4} and \ref{thn=lem-5}, we obtain a bijection $\phi$ between $A(n)$ and $B(n)$ as given by
\begin{equation}
\phi(\lambda)=\begin{cases}
\phi_1(\lambda),&\text{if }\lambda\in A_1(n);\\
\phi_2(\lambda),&\text{if }\lambda\in A_2(n);\\
\phi_3(\lambda),&\text{if }\lambda\in A_3(n);\\
\phi_4(\lambda),&\text{if }\lambda\in A_4(n);\\
\phi_5(\lambda),&\text{if }\lambda\in A_5(n).\\
\end{cases}
\end{equation}

 Now we are ready to  define the modified $pdt$-rank on partitions with overline designated summands.
 \begin{defi}\label{rmdtDef}
 Let $\lambda$ be a partition with overline designated summands. The modified $pdt$-rank of $\lambda$, denoted
 $r_{mdt}(\lambda)$,  is defined by
\begin{equation*}
r_{mdt}(\lambda)=\begin{cases}
1,&\text{if }\lambda\in A(n);\\
-1,&\text{if }\lambda\in B(n );\\
r_{dt}(\lambda),&\text{otherwise,}
\end{cases}
\end{equation*}
where $r_{dt}(\lambda)$ is the $pdt$-rank of $\lambda$.
\end{defi}

{\makeatletter\def\@captype{table}\makeatother\centering
\begin{tabular}{m{4cm}<{\centering}|m{4cm}<{\centering}|m{1.7cm}<{\centering}|m{1.7cm}<{\centering}}
$\lambda$       & $(\alpha,\beta)$ &  $r_{mdt}(\lambda)$  &$r_{mdt}(\lambda) \pmod 3$\\
\hline
$  \bar{5}'$  &  $(5,\ \emptyset)$
&  $1$        &  $1$\\
\hline
$  \bar{4}'+1'$     &  $(4,\ 1)$
&  $-1$             &  $2$\\
\hline
$4'+\bar{1}'$        &  $(1,\ 4)$
&$4$                 &  $1$\\
\hline
$\bar{3}'+2'$        &  $(3,\ 2)$
&  $2$               &  $2$\\
\hline
$3'+\bar{2}'$             &  $(2,\ 3)$
&  $3$                    &  $0$\\
\hline
$\bar{3}'+1'+1$      &  $(3,\ 1+1)$
&  $-2$              &  $1$\\
\hline
$3'+\bar{1}'+1$      &  $(1,\ 3+1)$
&  $0$               &  $0$\\
\hline
$\bar{3}'+1+1'$     &  $(3+1+1,\ \emptyset)$       &  $0$              &  $0$\\
\hline
$3'+1+\bar{1}'$     &  $(1+1,\ 3)$
&  $3$              &  $0$\\
\hline
$\bar{2}'+2+1'$     &  $(2,\ 2+1)$
&  $1$              &  $1$\\
\hline
$2'+2+\bar{1}'$     &  $(1,\ 2+2)$
&  $2$              &  $2$ \\
\hline
$2+\bar{2}'+1'$     &  $(2+2,\ 1)$
&  $-1$             &  $2$\\
\hline
$2+2'+\bar{1}'$     &  $(2+2+1,\emptyset )$       &  $0$              &  $0$\\
\hline
$\bar{2}'+1'+1+1$     &  $(2,\ 1+1+1)$
&  $-3$               &  $0$\\
\hline
$2'+\bar{1}'+1+1$     &  $(1,\ 2+1+1)$
&  $1$                &  $1$\\
\hline
$\bar{2}'+1+1'+1$     &  $(2+1+1,\ 1)$
&  $-1$               &  $2$\\
\hline
$\bar{2}'+1+1+1'$     &  $(2+1+1+1,\emptyset)$       &  $1$                &  $1$\\
\hline
$2'+1+1+\bar{1}'$     &  $(1+1+1,\ 2)$
&  $2$                &  $2$\\
\hline
$2'+1+\bar{1}'+1$     &  $(1+1,\ 2+1)$
&  $0$                &  $0$\\
\hline
$\bar{1}'+1+1+1+1$     &  $(1,\ 1+1+1+1)$       &  $-4$                &  $2$\\
\hline
$1+\bar{1}'+1+1+1$     &  $(1+1,\ 1+1+1)$       &  $-3$                &  $0$\\
\hline
$1+1+\bar{1}'+1+1$     &  $(1+1+1,\ 1+1)$       &  $-2$                &  $1$\\
\hline
$1+1+1+\bar{1}'+1$     &  $(1+1+1+1,\ 1)$       &  $-1$                &  $2$\\
\hline
$1+1+1+1+\bar{1}'$     &$(1+1+1+1+1,\emptyset)$       &  $1$                 &  $1$\\
\end{tabular}
 \caption{ The case for $n=5$ with modified $pdt$-rank $r_{mdt}(\lambda)$.}\label{pdt-r}
}

For example, for $n=5$, we have $PD_t(5)=24$. In Table \ref{pdt-r}, we
list the 24 partitions
of $5$ with overline designated summands, the corresponding pairs of
partitions along with the modified $pdt$-rank modulo $3$.
Let $N_{mdt}(i,t;n)$ denote the number of partitions of $n$ with overline designated summands with modified $pd$-rank congruent to $i \pmod{t}$.
It can be checked that
\[ N_{mdt}(0,3;5)=N_{mdt}(1,3;5)=N_{mdt}(2,3;5)=8.\]

 \vspace{.2cm}

\noindent{\bf Acknowledgments.}

The first author was supported by the Scientific Research Foundation of Nanjing Institute of Technology
(No. YKJ201627). The second author was supported by the National
Natural Science Foundation of China (No. 11801139) and the Natural Science Foundation of Jiangsu Province of China (No. BK20160855).
The third author was supported by the National
Natural Science Foundation of China (No. 11801119).

\end{document}